\documentclass[a4paper,11pt,reqno]{amsart}

\usepackage{amsmath,graphics}
\usepackage{amssymb}
\usepackage{amsfonts}
\usepackage{latexsym}
\usepackage{eucal}
\usepackage[dvips]{graphicx}

\theoremstyle{plain}
\newtheorem{thm}{Theorem}[section]
\newtheorem{propo}[thm]{Proposition}
\newtheorem{lem}[thm]{Lemma}

\theoremstyle{definition}

\renewcommand{\Re}{{\rm Re}}
\renewcommand{\Im}{{\rm Im}}
\newcommand{\R}{\mathbb{R}}
\newcommand{\C}{\mathbb{C}}

\newcommand{\N}{\mathbb{N}}

\renewcommand{\H}{\mathbb{H}}

\newcommand{\E}{\mathbb{E}}
\newcommand{\G}{\mathbf{G}}

\title[]{Random covers of compact surfaces and smooth linear spectral statistics}

\author[F.~Naud]{By Fr\'ed\'eric Naud }
\address{%
Fr\'ed\'eric Naud\\
Institut Math\'ematique de Jussieu\\
Universit\'e Pierre et Marie Curie, 4 place Jussieu, 75252 Paris Cedex 05\\
France.
}
\email{frederic.naud@imj-prg.fr}



\begin{document}
\bibliographystyle{plain}

\dedicatory{This paper is dedicated to the memory of Steve Zelditch, pioneer global analyst who introduced quantum ergodicity into the realm of pure mathematics.}
\maketitle

\begin{abstract}
We consider random $n$-covers $X_n$ of an arbitrary compact hyperbolic surface $X$. We show that in the large $n$ regime and small window limit, the variance of the smooth spectral statistics of the Laplacian $\Delta_\rho$ twisted by a unitary abelian character, obey the universal laws of GOE and GUE random matrices, depending on wether the character $\rho$ preserves or breaks the time reversal symmetry. 
We also prove a generalization for higher dimensional twists valued in compact linear groups.
These results confirm a conjecture of Berry \cite{Berry1,Berry2} and is a discrete analog of a recent work of Rudnick \cite{Rudnick} for the Weil-Petersson model of random surfaces.
\end{abstract}

\bigskip \tableofcontents

\section{Introduction and main results}
It is a long standing conjecture that the fine spectral statistics of classically chaotic quantum systems should follow, in the semi-classical regime, a universal behavior related to random matrices. In particular, Berry gave some theoretical evidence \cite{Berry1,Berry2} that time reversible systems, in a well chosen spectral window, should follow the spectral statistics of the Gaussian orthogonal ensemble (GOE), while non-reversible systems should follow the statistics of the Gaussian unitary ensemble (GUE). In the mathematical litterature, there are almost no rigourous results known, except the negative answer to Berry's conjecture in the case of arithmetic surfaces by Luo and Sarnak \cite{LuoSarnak} and the recent preprint of Rudnick \cite{Rudnick}. Rudnick's approach builds on the Weil-Petersson model of random hyperbolic surfaces and the celebrated integration formulas of Mirzakhani \cite{Mirza}  and most importantly Mirzhakani-Petri \cite{MirPet}.
In this paper, instead of averages over the whole moduli space, we will perform averages over the finite space of covers of degee $n$ and use the notion of random covers of compact hyperbolic surfaces 
defined by Magee, Puder and the author to study spectral gaps in \cite{MP1, MNP}.

Let $\H$ denote the hyperbolic plane endowed with its metric of constant curvature $-1$, and let $\Gamma \subset \mathrm{PSL}_2(\R)$ be a co-compact Fuchsian group so that $X=\Gamma \backslash \H$ is a compact hyperbolic surface. Let $\phi_n:\Gamma \rightarrow \mathcal{S}_n$ be a group homomorphism, where $\mathcal{S}_n$ is  the symmetric group of permutations of $[n]:=\{1,\ldots,n\}$. The discrete group $\Gamma$ acts on $\H\times [n]$ by
$$\gamma.(z,j):=(\gamma(z), \phi_n(\gamma)(j)). $$
The resulting quotient $X_n:=\Gamma \backslash \H\times [n]$ is then a finite cover of degree $n$ of $X$, possibly not connected. By considering the (finite) space of all homomorphism $\phi_n:\Gamma\rightarrow \mathcal{S}_n$, endowed with the uniform probability measure, we obtain a notion of Random covering surfaces of degree $n$, $X_n\rightarrow X$. The expectation of any random variable $Y$ with respect to this uniform measure is denoted by $\mathbb{E}_n(Y)$. We point out that we can also view (up to isometry) the random cover $X_n$ as 
$$X_n=\bigsqcup_{k=1}^p \Gamma_k\backslash \H,$$
where each $\Gamma_k$ is (in general non-normal) subgroup of $\Gamma$ given by $\Gamma_k=\mathrm{Stab}_\Gamma(i_k)=\{ \gamma \in \Gamma\ :\ \phi_n(\gamma)(i_k)=i_k \}$ where $i_1,\ldots,i_k \in [n]$ are representatives of the orbits of $\Gamma$ (acting on $[n]$ via $\phi_n$). In general, the cover $X_n$ is not connected, but it follows directly from Liebeck and Shalev \cite{LS1} that the probability that this cover is connected (i.e. $\Gamma$ acts transitively on $[n]$ via $\phi_n$) tends to $1$ as $n$ goes to infinity.

\bigskip In order to allow us to observe different statistical behaviors, we also fix a complex unitary representation, globally defined on $\Gamma$,  
$$\rho:\Gamma \rightarrow GL(V_\rho),$$
where $V_\rho$ is a finite-dimensional complex vector space.

Each surface $X_n$ comes with a self-adjoint Laplace operator $\Delta_{n,\rho}$, the hyperbolic Laplacian on $X_n$ twisted by  $\rho$. The spectrum of $\Delta_{n,\rho}$ is denoted  (with multiplicity) by
$$\left \{ \lambda_{j,n}=\frac{1}{4}+(r_{j,n})^2,\  j\in \N \right \},$$
where $r_{j,n}\in \R \cup i\R$.

Let $\psi$ be a real valued {\it  even} smooth test function on $\R$ whose Fourier transform $\widehat{\psi}$ is compactly supported in $\R$. The Fourier transform is in our case given by
$$\widehat{\psi}(\xi)=\frac{1}{2\pi}\int_{-\infty}^{+\infty} \psi(x)e^{-ix\xi} dx.$$

Note that $\psi$ is therefore an entire holomorphic function on $\C$. We take $\alpha, L>0$ and consider the smooth counting function
$$N_{n}(L):=\sum_{j=0}^\infty \left \{ \psi(L(r_{j,n}-\alpha))+\psi(L(r_{j,n}+\alpha) ) \right \}.$$
We will fix $\alpha$ (the energy level) and let $L^{-1}$ (the band width) go to zero and $n$ go to infinity.
It is easy to see that for fixed $L$, we have as $n$ goes to infinity,
$$ \E_n(N_n(L))\sim C_\alpha\frac{(g-1)\mathrm{dim}(V_\rho)n}{L}\int_\R \psi(r)dr ,$$
where $g$ is the genus of the base surface $X$, and $C_\alpha=2\alpha \tanh(\pi \alpha)$. One of the simplest quantity to investigate is the so-called "number variance" which measures
"spectral rigidity". We will therefore focus on the probabilistic variance
$$\mathbb{V}_n(L):=\mathbb{E}_n \left (\left \vert N_n(L)-\mathbb{E}_n(N_n(L)) \right \vert^2\right).$$
Our main result is the following.
\begin{thm} 
\label{main}
Fix $X=\Gamma \backslash \H$ and a unitary abelian character $\rho:\Gamma \rightarrow \C$ as above. We fix an energy level $\alpha>0$, then we have the double limit
$$\lim_{L\rightarrow \infty} \lim_{n\rightarrow \infty} \mathbb{V}_n(L)=
\left \{ \Sigma^2_{\mathrm{GOE}}(\psi)\ \mathrm{if}\  \rho^2=1\ \atop \Sigma^2_{\mathrm{GUE}}(\psi)\ \mathrm{if}\ \rho^2 \neq 1. \right.$$
The quantity $\Sigma^2_{\mathrm{GOE}}(\psi)$ is the "smoothed" number variance of random matrices for the GOE model in the large dimension limit and given by
$$\Sigma^2_{\mathrm{GOE}}(\varphi)=2\int_{\R} \vert x \vert (\widehat{\psi}(x))^2dx,$$
while $\Sigma^2_{\mathrm{GUE}}(\psi)=\frac{1}{2}\Sigma^2_{\mathrm{GOE}}(\psi)$.
\label{main1}

\end{thm}
This result is in accordance with Berry's conjecture: the hypothesis $\rho^2\neq 1$ corresponds to breaking the time reversal symmetry i.e. the character $\chi$ is sensitive to the orientation of closed geodesics.
The computation of number variances in the model cases of GUE, GOE goes back to Dyson and Mehta in \cite{DM1}, see Mehta \cite{MehtaBook}, Appendix A.39 and also chapter 16.  

\bigskip
What about higher dimensional twists and Laplacians acting on associated flat vector bundles? We denote by $U(N)$ the unitary group, $SU(N)$ the special unitary group, $Sp(N)$ the compact symplectic group, $SO(N)$ the special orthogonal group. Note that $SU(N), SO(N)$ are Lie subgroups of $U(N)$ while for the compact symplectic group
$$Sp(N)=Sp(2N,\C)\cap U(2N)$$ is a Lie subgroup of $U(2N)$. In the text below $\G$ will denote a compact connected real Lie subgroup of $U(N)$ for some $N$. Notice that by a celebrated
theorem of Tannaka, see Chevalley \cite{Chevalley} , this is automatically a real algebraic group.
We show the following .
\begin{thm}
\label{main2}
Let $\G$ be as above, and assume it is not conjugated to a direct product of Lie groups.
Fix $X=\Gamma \backslash \H$ and a unitary representation $\rho:\Gamma \rightarrow \G$ whose image is Zariski dense in $\G$. We fix an energy level $\alpha>0$, then we have the double limit
$$\lim_{L\rightarrow \infty} \lim_{n\rightarrow \infty} \mathbb{V}_n(L)=\left\{ \Sigma_{\mathrm{GOE}}^2(\psi)\ \mathrm{if}\ g\mapsto \mathrm{Tr}(g)\in \R, \atop \Sigma^2_{\mathrm{GUE}}(\psi)\ \mathrm{otherwise}. \right.$$
In particular, we have the following behavior for classical compact Lie groups.
\begin{itemize}
\item For $\G=U(N)$, or $\G=SU(N)$ with $N\geq 3$, we are in the GUE case.
\item For $\G=SO(N),Sp(N)$, we are in the GOE case.
\end{itemize}
\end{thm}
In a nutshell, if the trace map on $\G$ is {\it real-valued} then we get GOE statistics while if the trace map takes complex values, we have GUE statistics. Representation of surface groups
with Zariski dense images in the above cases do exist under some mild dimensional assumptions, see for example in the paper of Kim and Pansu, \cite{KP}  and also Kishore \cite{kishore}.
More details can be found at the end 
of the paper, $\S 5$.

\bigskip \noindent
While our results have the same flavor as in \cite{Rudnick}, here are the main differences:
\begin{itemize}
\item We take finite averages over (increasingly bigger) sets of surfaces. In particular we cannot rely as much on oscillatory integral techniques and integration by parts as in the smooth model of Weil-Petersson.
\item Our result is at fixed frequency $\alpha$ while Rudnick in \cite{Rudnick} lets both $L=O(\log\alpha)),\alpha$ go to $\infty$ . However, it is not really necessary at the technical level in \cite{Rudnick} to let the energy level $\alpha$ go to infinity in order to observe the universal behaviour for the variance so both result are actually the same.
\item Our model also allows to study twisted Laplacians and prove the conjectured effect of time reversal symmetry on the spectral statistics. It would be interesting to know to what extent Mirzakhani's formalism needs to be modified to allow unitary twists in the analysis of Rudnick.
\item If $\Gamma$ is an arithmetic groups, then all the random covers $X_n$ are also arithmetic surfaces. It was shown in \cite{LuoSarnak}, see also \cite{BS}, that {\it in the semi-classical} regime, the number variance exhibits some poissonian behaviour, due to the high multiplicities in lengths of closed geodesics. Our result shows that in the large $n$ regime we still get GOE/GUE statistics despite arithmeticity.
\end{itemize}

The paper is organized as follows. In $\S2$ we explain how the spectrum of $\Delta_{n,\rho}$ can be studied via an $n\times \mathrm{dim}(V_\rho)$-dimensional unitary twist of the Laplacian on the base surface $X$, mostly by using the induction formula from representation theory.
We then state the ad-hoc Selberg trace formula that we will later use in the proof. In $\S 3$ we review the asymptotic statistics (in the large $n$ regime) of the number of fixed points of "random" permutation $\phi_n(\gamma)$: when $\gamma$ is primitive, this random variable behaves essentially like a Poisson variable with parameter $1$. When $\gamma_1\neq \gamma_2$ are in two different primitive conjugacy class (and are not inverses of each other), the associated number of fixed points of $\phi_n(\gamma_1)$ and $\phi_n(\gamma_2)$ are asymptotically independent. This will play a key role in the proof of Theorem \ref{main}. We point out that a similar behaviour occur in the smooth Weil-Peterssen model where the core argument in \cite{MirPet} shows independency and Poissonian behaviour for lengths of simple closed geodesics.

In $\S 4$ we provide some asymptotic results for the (weighted) counting function for the closed geodesics on $X$, which we derive from some standard
facts on twisted Laplaclans for unitary representations and the trace formula. Finally, $\S 5$ is devoted to the proof of Theorems \ref{main} and \ref{main2} which in turn is broken into several steps. 

\bigskip \noindent
{\bf Acknowledgements}. The author wants to thank Michael Magee, St\'ephane Nonnenmacher, Doron Puder and Zeev Rudnick for several discussions around this work.
\section{Representations, twisted laplace spectrum and trace formula}
Let $\rho:\Gamma \rightarrow  U(V_\rho)$ be a {\it unitary representation} of $\Gamma$, with finite dimensional complex representation space $V_\rho$, and inner hermitian product $\langle .,. \rangle_\rho$. The twisted Laplacian $\Delta_\rho$ is the Laplacian acting on 
smooth functions $$\psi:\H\rightarrow V_\rho$$ satisfying the equivariance property
$$\forall\ \gamma \in \Gamma,\ \psi(\gamma z)=\rho(\gamma)\psi(z).$$
We endow this space with an $L^2$-norm given by
$$\Vert \psi\Vert_\rho^2:=\int_{\mathcal F} \Vert \psi(z) \Vert_\rho^2 d\mathrm{Vol}(z),$$
where $\mathcal{F}\subset \H$ is a compact fundamental domain for $\Gamma$ and $\mathrm{Vol}$ is hyperbolic volume.
This Laplacian $\Delta_\rho$ has a unique self-adjoint extension with discrete spectrum.  Alternatively, one can view $\Delta_\rho$ as a self-adjoint Bochner-Laplace operator associated to a metric connection on the flat vector Bundle
$E_\rho\rightarrow X$ associated to the representation $\rho$ of the fundamental group $\Gamma$. Spectral theory and trace formulas for such Laplacians are well documented in \cite{Hejhal1,Venkov}.

 We have the following fact.
\begin{propo}
  There exists an $n\times \mathrm{dim}(V_\rho)$-dimensional complex representation 
  $$\varrho_n:\Gamma \rightarrow U(V_n)$$ such that the spectrum of $\Delta_{n,\rho}$ on the random surface $X_n$ coincides (with same multiplicity) with the spectrum of the twisted Laplacian 
  $\Delta_{\varrho_n}$. In addition, we have the character formula for all $\gamma \in \Gamma$, 
  $$\mathrm{Tr}(\varrho_n(\gamma))=\mathrm{Tr}(\rho(\gamma)) \vert \mathrm{Fix}(\phi_n(\gamma))\vert,$$
  where $\mathrm{Fix}(\phi_n(\gamma))$ denotes the set of fixed points in $[n]$ of the permutation $\phi_n(\gamma)$.
\end{propo}
{\it Proof.} We set $d=\mathrm{dim}(V_\rho)$. We recall that we have 
$$X_n=\bigsqcup_{k=1}^p \Gamma_k\backslash \H,$$
where $\Gamma_k=\mathrm{Stab}_\Gamma(i_k)=\{ \gamma \in \Gamma\ :\ \phi_n(\gamma)(i_k)=i_k \}$ where $i_1,\ldots,i_k \in [n]$ are representatives of the orbits of $\Gamma$ (acting on $[n]$ via $\phi_n$). Working on each connected components of $X_n$, it is enough to assume that $X_n$ is actually connected, i.e. $\Gamma$ acts transitively on $[n]$. Let us set
$$\Gamma_1:=\mathrm{Stab}_\Gamma(1),$$
and we write the right cosets decomposition of $\Gamma$ as
$$\Gamma=\Gamma_1 g_1 \cup \Gamma_1 g_2 \cup \ldots \cup \Gamma_1 g_n,$$
where $g_1=Id$ and we have ordered the cosets such that $\phi_n(g_i)(i)=1$ for all $i=1,\ldots,n$. It is a well known fact (see for example Venkov \cite{Venkov}, page 51) that the Laplacian $\Delta_\rho$ on $X_n=\Gamma_1\backslash \H$ is unitarily conjugated to the Laplacian on $E_{\varrho_n}$,
where $\varrho_n$ is the {\it induced representation} of $\rho$ on $\Gamma_1$ to $\Gamma$. More precisely, $\varrho_n(\gamma)$ acts on $V_n=\C^{d\times n}$ via the block-matrix 
$$\mathrm{Mat}_{\C^{d\times n}}(\varrho_n(\gamma))=\left [   \widetilde{\rho}(g_i\gamma g_j^{-1}) \right ]_{1\leq i,j\leq n},$$
where $\widetilde{\rho}(g)=\rho(g)$ if $g\in \Gamma_1$ and $\widetilde{\rho}(g)=0$ elsewhere.
On the diagonal part of this matrix, observe that if $\phi_n(\gamma)(i)\neq i$ then $\widetilde{\rho}(g_i\gamma g_i^{-1})=0$ while
$$\widetilde{\rho}(g_i\gamma g_i^{-1})=\rho(g_i \gamma g_i^{-1})=\rho(g_i)\rho(\gamma)\rho(g_i)^{-1}$$ if $\phi_n(\gamma)(i)=i$. Computing traces yields
$$\mathrm{Tr}(\varrho_n(\gamma))=\sum_{i\ :\ \phi_n(\gamma)(i)=i} \mathrm{Tr}(\rho(g_i)\rho(\gamma)\rho(g_i)^{-1})=\mathrm{Tr}(\rho(\gamma)) \vert \mathrm{Fix}(\phi_n(\gamma))\vert,$$
and the proof is done. $\square$

\bigskip We will use twisted versions of Selberg's trace formula which is as follows, see for example in \cite{Hejhal1, Venkov}. Let $\varrho:\Gamma \rightarrow U(V_\varrho)$ be a general unitary, finite dimensional, representation of $\Gamma$. 
For all eigenvalue $\lambda_j(\varrho)$ of the Laplacian $\Delta_\varrho$, we fix $r_j(\rho)$ such that $\lambda_j(\varrho)=1/4+(r_j(\varrho))^2$. Let $h$ be an even function whose Fourier transform
$\widehat{h}$ is compactly supported and smooth, then we have the celebrated identity:
\begin{thm}(Twisted trace formula)
\label{trace1}
$$\sum_{j=0}^\infty h(r_{j}(\varrho))=(g-1)\mathrm{dim}(V_\varrho)\int_{-\infty}^\infty h(r)r\tanh(\pi r)dr$$
$$+\sum_{\gamma \in \mathcal{P}} \sum_{k\geq 1}  \mathrm{Tr}(\varrho(\gamma^k))  \frac{\ell(\gamma)\widehat{h}(k\ell(\gamma))}{2\sinh(k\ell(\gamma)/2)}.$$
In the above formula, $g$ denotes the genus of the base surface $X$, $\mathcal{P}$ is the set of primitive conjugacy classes in $\Gamma$ (different from identity) and if $\gamma \in \mathcal{P}$,
$\ell(\gamma)$ denotes the length of the associated closed geodesic on $X$. The sum over $k$ runs over all positive integers. 
\end{thm}

\section{Character statistics in the large $n$ regime}
In this section we review some key facts on the integer valued random variables $\vert\mathrm{Fix}(\phi_n(\gamma))\vert$, for $\gamma \in \Gamma$ and where $\phi_n:\Gamma \rightarrow \mathcal{S}_n$ is a random homomorphism. To simplify further notations, we will denote simply $\vert \mathrm{Fix}(\phi_n(\gamma))\vert$ by $F_n(\gamma)$. Our analysis will based on the following facts. We recall that $\mathcal{P}$ denotes the set of primitive conjugacy classes (different from identity) in $\Gamma$.
\begin{propo}
\label{Stat1} Under the above notations, the following holds.
\begin{enumerate}
\item For all $\gamma \in \mathcal{P}$ and $k\geq 1$, we have as $n$ goes to $\infty$,
$$\E_n( F_n(\gamma^k))=d(k)+O\left( \frac{1}{n} \right),$$
where $d(k)$ denotes the number of divisors of $k$.
\item For all $\gamma_1,\gamma_2 \in \mathcal{P}$ such that $\gamma_1\not \in \{\gamma_2,\gamma_2^{-1} \}$, for all $k_1,k_2 \geq 1$, we have as $n\rightarrow \infty$,
$$\E_n( F_n(\gamma_1^{k_1}) F_n(\gamma_2^{k_2}))=\E_n( F_n(\gamma_1^{k_1}))\E_n( F_n(\gamma_2^{k_2}))+O\left( \frac{1}{n} \right).$$
\item For all $k_1,k_2\geq 1$ and $\gamma \in \mathcal{P}$, the limit 
$$\lim_{n\rightarrow +\infty} \E_n\left [  \left(F_n(\gamma^{k_1})-\E_n(F_n(\gamma^{k_1}))\right)\left(F_n(\gamma^{k_2})-\E_n(F_n(\gamma^{k_2}))\right)\right]$$
$$:=\mathcal{V}(k_1,k_2)$$
exists. If $k_1=k_2=1$ then $\mathcal{V}(1,1)=1$. Moreover  we have the formula
$$\mathcal{V}(k_1,k_2)=\sum_{d\vert k_1\& k_2}d.$$
\end{enumerate}
\end{propo}
{\it Proof}. From \cite{PZ,MP1}, see for example in \cite{PZ} Theorem 1.4, we have directly $(1)$.  From \cite{PZ}, Theorem 1.13, we know that if $\gamma_1,\gamma_2$ do not belong to the same cyclic subgroup, then the random variables
$F_n(\gamma_1),F_n(\gamma_2)$ are asymptotically independent which implies $(2)$. Denote by $C_{n,d}(\gamma)$ the number of $d$-cycles of the random permutation $\phi_n(\gamma)$,
where $\gamma$ is primitive.  As usual $Z_\lambda$ will denote a Poisson variable with parameter $\lambda>0$ whose law is given by
$$\mathbb{P}(Z_\lambda=q)=\frac{e^{-\lambda} \lambda^q}{q!},$$
and satisfies $\E(Z_\lambda)=\lambda$, $\E(Z_\lambda^2)=\lambda^2+\lambda$.
We know from \cite{PZ,MP1}, Theorem 1.14 that as $n$ goes to $\infty$, $C_{n,d}(\gamma)$ converges in distribution to a Poisson variable $Z_{1/d}$. We point out that this result was 
already known for free groups since Nica \cite{Nica}. 

Moreover, if $d_1\neq d_2$ then the variables
$C_{n,d_1}(\gamma),C_{n,d_2}(\gamma)$ are asymptotically independent. Using the formula (for $\gamma$ primitive) 
$$F_n(\gamma^k)=\sum_{d\vert k} d C_{n,d}(\gamma),$$
we get that
$$\E_n( F_n(\gamma^{k_1})F_n(\gamma^{k_2}))=\sum_{d_1\vert k_1}\sum_{d_2 \vert k_2}d_1d_2 \E_n(C_{n,d_1}(\gamma)C_{n,d_2}(\gamma)).$$
By the Above facts, if $d_1\neq d_2$ then 
$$\lim_{n\rightarrow +\infty} \E_n(C_{n,d_1}(\gamma)C_{n,d_2}(\gamma))=\E(Z_{1/d_1})\E(Z_{1/d_2})=\frac{1}{d_1d_2}.$$
On the other hand, if $d_1=d_2$ then
$$ \lim_{n\rightarrow +\infty} \E_n(C_{n,d_1}(\gamma)C_{n,d_2}(\gamma))=\E(Z_{1/d_1}^2)=\frac{1}{d_1}+\frac{1}{d_1^2}.$$
Using the fact that 
$$\lim_{n\rightarrow +\infty} \E_n(F_n(\gamma^k))=d(k)=\sum_{d\vert k} 1=\sum_{d\vert k} d \lim_{n\rightarrow \infty} \E_n(C_{n,d}(\gamma)),$$
we get that
$$\lim_{n\rightarrow +\infty} \E_n\left [  \left(F_n(\gamma^{k_1})-\E_n(F_n(\gamma^{k_1}))\right)\left(F_n(\gamma^{k_2})-\E_n(F_n(\gamma^{k_2}))\right)\right]$$
$$=\sum_{d_1\vert k_1} \sum_{d_2 \vert k_2} d_1d_2\left (  \lim_{n\rightarrow +\infty} \E_n(C_{n,d_1}(\gamma)C_{n,d_2}(\gamma))-\frac{1}{d_1d_2}  \right).$$
$$=\sum_{d\vert k_1 \& k_2} d,$$
and the proof is done. $\square$

\section{Bounds on closed geodesics and Chebotarev counting results}
We will use throughout this paper the following standard summation by parts formula: this is standard in any introductory class in analytic number theory.
\begin{propo} 
Let $0<\ell_0<\ell_1<\ldots<\ell_n$ be an increasing sequence of positive real numbers with $\lim_{n\rightarrow \infty} \ell_n=+\infty$, and let 
$(c_n)_{n\in \N}$ be an arbitrary sequence of complex numbers. Let $f:(0,+\infty)\rightarrow \C$ be a $C^1$-function. We have for all $x\geq 0$,
$$\sum_{\ell_n\leq x} c_nf(\ell_n)=A(x)f(x)-\int_0^x A(t)f'(t)dt, $$
where 
$$A(x)=\left \{ \sum_{\ell_n\leq x}c_n\ \mathrm{if}\ x\geq \ell_0 \atop 0\ \mathrm{otherwise.} \right.$$
\end{propo}
 It will be sometimes convenient to use the formal notation
$$\int_0^x f(t)dA(t):=\sum_{\ell_n\leq x} c_nf(\ell_n),$$
which in the special case when $(c_n)$ is positive, coincides with an actual Stieltjes integral.
\subsection{Some counting bounds}
Let $N(T)$ denote the counting function for closed geodesics on $X=\Gamma \backslash X$ whose length is smaller than $T$ i.e.
$$N(T):=\sum_{k\ell(\gamma))\leq T} 1, $$
where the sum runs over all integers $k\geq 1$ and all $\gamma \in \mathcal{P}$. One can also consider the primitive counting function 
$$N^0(T):=\sum_{\ell(\gamma))\leq T} 1,$$
where the sum runs only on primitive conjugacy classes.
It has been known since Huber \cite{Huber} that as $T\rightarrow +\infty$,
$$N^0(T)=\frac{e^T}{T}(1+o(1)),$$
which is commonly referred as "the prime geodesic theorem".
 We point out that since we have
$$N(T)=N^0(T)+\sum_{k\geq 2}N^0(T/k), $$
we actually have as $T\rightarrow +\infty$,
$$N(T)=N^0(T)+O(Te^{T/2}).$$
Later \cite{Sarnak1,Hejhal1} more precise remainders where proved: for all $\epsilon>0$, we have the following asymptotic formula as $T$ goes to infinity
\begin{equation}
\label{prime0}
 N(T)=\mathrm{Li}\left (e^{ T}\right)+O\left(e^{\alpha_0(\epsilon) T} \right),
 \end{equation}
 where $\alpha_0(\epsilon)=\max\{ 3/4+\epsilon, s_1(X)\}$, where $0<s_1(X)<1$ is related to the first non-trivial eigenvalue of the Laplacian on $X$ by $s_1(1-s_1)=\lambda_1(X)$.
 From the fact that
 $$\mathrm{Li}(x):=\int_2^x \frac{dt}{\log t},$$
 it is immediate to check that as $T\rightarrow +\infty$,
 \begin{equation}
 \label{prime1}
 N^0(T)=\frac{e^T}{T}\left ( 1+O\left( \frac{1}{T}\right) \right),
 \end{equation}
 which is sometimes sufficient. Given a unitary representation 
 $$\rho:\Gamma\rightarrow GL(V_\rho), $$
 we set $\chi_0(\gamma):=\mathrm{Tr}(\rho(\gamma))$. We then consider the following counting function
 $$N_{\chi_0}(T):=\sum_{k\ell(\gamma))\leq T} \chi_0(\gamma^k). $$
 We will show the following fact.
 \begin{propo}
 \label{prime2} Assume that $\rho$ does not contain the trivial representation, then there exists $\beta_0=\beta_0(\Gamma,\chi_0)<1$ such that as $T\rightarrow +\infty$,
 $$\vert N_{\chi_0}(T)\vert=O\left(e^{\beta_0 T} \right).$$
 \end{propo}
 \noindent {\it Proof}. The proof follows from standard facts on the spectrum of the twisted Laplacian $\Delta_{\rho}$ and will certainly not surprise the experts. We nevertheless include a concise proof for completeness. 
 Let $0\leq \lambda_0(\rho)\leq \lambda_1(\rho)\leq \ldots$ denote the spectrum of $\Delta_{\rho}$. From the analysis of Phillips-Sarnak \cite{PS1} in the scalar case or from Sunada \cite{Sunada}  in the general vector case, we know that if $\rho$ does not contain the trivial representation, then we have $\lambda_0(\rho)>0$. To be more precise, it is actually shown in \cite{Sunada} that if
 $\mathcal{A}$ is a finite set of generators of $\Gamma$ then there exists $C(\mathcal{A},\Gamma)>0$ such that
 $$C(\mathcal{A},\Gamma)\delta_{\mathcal A}(\rho,{\mathbf 1})^2\leq \lambda_0(\rho), $$
where $\delta_{\mathcal A}(\rho,{\mathbf 1})$ is Khazdan's distance to the trivial representation $\bf 1$. Khazdan's distance is defined by
$$\delta_{\mathcal A}(\rho,{\mathbf 1})=\inf_{v\in V_\rho \atop \Vert v\Vert=1} \max_{\gamma \in \mathcal{A}} \Vert \rho(\gamma)v-v\Vert. $$
 It is straightworward to check that if $\delta_{\mathcal A}(\rho,{\mathbf 1})=0$, then one can find $v\neq 0 \in V_\rho$ such that for all $\gamma \in \Gamma$, $\rho(\gamma)v=v$ i.e. $\rho$ contains the trivial representation.

 \bigskip
 The end of the proof now follows from an application of the trace formula. First observe that we have
 $$ N_{\chi_0}(T)=\sum_{k\ell(\gamma))\leq T} \chi_0(\gamma^k)=\sum_{\gamma \in \mathcal{P} \atop \ell(\gamma)\leq T} \chi_0(\gamma) +O\left( Te^{T/2}\right),$$
 just by noticing that by the prime orbit theorem we have
 $$ \sum_{k\ell(\gamma))\leq T \atop k\geq 2} 1=O\left( Te^{T/2}\right).$$
 Therefore it is enough to prove the desired result for 
 $$N_{\chi_0}^0(T):=\sum_{\ell(\gamma))\leq T} \chi_0(\gamma).$$
 
 Let $\ell_0>0$ denote the shortest length of closed geodesics on $X$. Pick $\ell_0<T_0<T_0+T_1$, and let $0\leq q(x)\leq 1$ be an even smooth compactly supported function
 on $\R$ such that $q(x)=0$ for all $x\in [-\ell_0/2,+\ell_0/2]$ and $q(x)=1$ for all $\ell_0\leq \vert x\vert \leq T_0$ while being supported inside $\{ \vert x\vert \leq T_0+T_1 \}$. 
 In the latter we will choose $T_1$ as a function of $T_0$ satisfying $T_1=e^{-\alpha T_0}$, for some $\alpha>0$, while all the bounds will be relevant for $T_0\rightarrow +\infty$. 
 Moreover, it is possible to choose $q$ such that for all $k$,
 $$\sup_\R \vert q^{(k)}\vert \leq C_k T_1^{-k}.$$
 Consider now the test function
 $$\varphi_{T_0,T_1}(x):=\frac{2\cosh(x/2)}{\vert x\vert}q(x),$$
and set 
$$h_{T_0,T_1}(u):=\int_\R e^{iux} \varphi_{T_0,T_1}(x)dx.$$ 
Applying the trace formula from Theorem \ref{trace1},  a direct bound shows that the finite contribution from the low eigenvalues below $1/4$
satisfies 
$$\sum_{j\ :\  r_j \in i\R} \vert h_{T_0,T_1}(r_j) \vert=O \left( T_0 e^{(1/2+s_0)T_0}    \right),$$
where $$s_0=\left \{ \sqrt{1/4-\lambda_0}\ \mathrm{if}\ \lambda_0\leq 1/4, \atop 0\  \mathrm{otherwise}. \right.$$
Note that since $\lambda_0>0$, we have $s_0<1/2$. Integrating by parts 3 times yields
$$h_{T_0,T_1}(u)=O\left( \frac{e^{T_0/2}T_1^{-3}}{1+\vert u \vert^3} \right),$$
which gives by a crude bound
$$\int_{-\infty}^\infty h_{T_0,T_1}(r)r\tanh(\pi r)dr=O\left ( e^{T_0/2}T_1^{-3}\right).$$
We recall that by Weyl's law we have the bound
$$\#\{ j\ :\ \vert r_j\vert \leq R\}=O(R^2),$$
so that a standard summation by parts yields also the crude bound
$$\sum_{j\ :\ r_j\in \R} \vert h_{T_0,T_1}(r_j)\vert=O\left ( e^{T_0/2}T_1^{-3}  \right).$$
We now observe that
$$ \sum_{\gamma \in \mathcal{P}} \sum_{k\geq 1}  \chi_0(\gamma^k)  \frac{\ell(\gamma)\varphi_{T_0,T_1}(k\ell(\gamma))}{2\sinh(k\ell(\gamma)/2)}=
\sum_{\gamma \in \mathcal{P}}   \chi_0(\gamma^k)(\tanh(\ell(\gamma)/2))^{-1}q(\ell(\gamma))+O(T_0 e^{T_0/2}).$$
By the  the fact that 
$$\frac{1}{\tanh(x)}=1+O(e^{-2x}),$$
we have 
$$\sum_{\gamma \in \mathcal{P}}   \chi_0(\gamma^k)(\tanh(\ell(\gamma)/2))^{-1}q(\ell(\gamma))=\sum_{\gamma \in \mathcal{P}}   \chi_0(\gamma^k)q(\ell(\gamma)) 
+O\left( \sum_{\gamma} e^{-\ell(\gamma)}q(\ell(\gamma)) \right),$$
and the latter term is easily shown to be $O(1)$ by summation by parts. By the very definition of $q$ we have obtained that
$$\sum_{\gamma \in \mathcal{P}} \sum_{k\geq 1}  \chi_0(\gamma^k)  \frac{\ell(\gamma)\varphi_{T_0,T_1}(k\ell(\gamma))}{2\sinh(k\ell(\gamma)/2)}$$
$$=N^0_{\chi_0}(T_0)+
O(T_0 e^{T_0/2})+\sum_{T_0\leq \ell(\gamma)\leq T_0+T_1} \chi_0(\gamma^k)q(\ell(\gamma)). $$
By the prime orbit theorem with remainder (\ref{prime0}),
$$\left \vert \sum_{T_0\leq \ell(\gamma)\leq T_0+T_1} \chi_0(\gamma^k)q(\ell(\gamma)) \right \vert\leq O\left (N^0(T_0+T_1)-N^0(T_0)\right)=O(e^{T_0}T_1)+O(e^{\alpha_0 T_0}).$$
Gathering all these estimates, we have thus
$$ N^0_{\chi_0}(T_0)=O(T_0 e^{T_0/2})+O(e^{T_0}T_1)+O(e^{\alpha_0 T_0})+O\left ( e^{T_0/2}T_1^{-3}\right)+O \left( T_0 e^{(1/2+s_0)T_0}    \right),$$
where $\alpha_0<1$, and $s_0<1/2$. Setting $\alpha=\frac{1}{8}$ gives
$$N^0_{\chi_0}(T_0)=O(T_0 e^{T_0/2})+O(e^{\alpha_0 T_0})+O \left( T_0 e^{(1/2+s_0)T_0} \right) +O\left(e^{\frac{7}{8}T_0} \right),$$
and the proof is done. $\square$

\bigskip
Notice that we have not attempted at all to optimize the error term since any exponent $\beta_0<1$ will be enough for our application. We will now apply this result to prove the following
asymptotic formula.
\begin{propo}(Chebotarev type asymptotic)
\label{equi2} Assume that $\rho:\Gamma\rightarrow \mathbf{G}$ is a unitary representation where $\mathbf{G}$ is a compact algebraic Lie subgroup of $GL_N(\C)$. Assume that the image of $\Gamma$ by $\rho$ is Zariski dense in $\mathbf{G}$. Let $f(g)=\phi(\mathrm{Tr}(g))$ be a polynomial function of the trace. Then as $x\rightarrow +\infty$ we have
$$\sum_{\ell(\gamma)\leq x} f(\rho(\gamma))= \left (\int_{\G}f(g)dg\right)\mathrm{Li}(e^x)+O(e^{\alpha_0 x}),$$
for some $\alpha_0<1$ and where $dg$ denotes normalized Haar measure on $\G$.
\end{propo}
\noindent {\it Proof}. By Peter-Weyl theorem, since $f(g)$ is a smooth central function we know that we can expand it as
$$f(g)=\sum_{\lambda \in \widehat{\G}} c_\lambda \mathrm{Tr}(\lambda(g)),$$
where $c_\lambda \in \C$ the sum runs over all irreducible representations $\lambda$ of $\G$. This sum is actually finite: indeed by the polynomial hypothesis we have the finite sum
$$f(g)=\sum_{k,l} a_{k,l} (\mathrm{Tr}(g))^k (\overline{\mathrm{Tr}(g)})^l,$$
where each $a_{k,l}\in \C$. Let $\mathcal{I}:\G\rightarrow GL(\C^N)$ denote the faithfull representation of $\G$ acting as usual on $\C^N$ by left matrix product with column vectors. Then we have
$$(\mathrm{Tr}(g))^k (\overline{\mathrm{Tr}(g)})^l=\mathrm{Tr}\left ( \left(\mathcal{I}^{\otimes k} \otimes     \widehat{\mathcal{I}}^{\otimes l}\right)(g)   \right ),$$
where $\widehat{\mathcal{I}}$ is the contragredient representation of $\mathcal{I}$, see for example \cite{Bump} Proposition 2.6. The finite-dimensional representation
$$\mathcal{I}^{\otimes k} \otimes     \widehat{\mathcal{I}}^{\otimes l}$$
can in turn be written as a finite direct sum of irreducible representations of $\G$, hence $(\mathrm{Tr}(g))^k (\overline{\mathrm{Tr}(g)})^l$ can be expressed as a finite linear combination
of characters $\mathrm{Tr}(\lambda(g))$, and so does $f(g)$. The uniqueness of Peter-Weyl decomposition then implies that only finitely many coefficients $c_\lambda$ must be non-vanishing, notice in addition that we have 
$$c_{\mathrm{trivial}}=\int_\G f(g)dg.$$
We then write
$$\sum_{\ell(\gamma)\leq x} f(\rho(\gamma))=\sum_{\lambda \in \widehat{\G}} c_\lambda \sum_{\ell(\gamma)\leq x}  \mathrm{Tr}(\lambda(\rho(\gamma)))$$
$$= \left (\int_{\G}f(g)dg\right)\sum_{\ell(\gamma)\leq x} 1+ \sum_{\lambda \in \widehat{\G}\atop \lambda\ \mathrm{non\ trivial}} c_\lambda \sum_{\ell(\gamma)\leq x}  \mathrm{Tr}(\lambda(\rho(\gamma))).$$
Assume now that if $\lambda$ is irreducible and assume that $\lambda \circ \rho$ contains the trivial representation. This means that there exists a non vanishing vector $v \in V_\lambda$
such that for all $\gamma \in \Gamma$, 
$$\lambda \circ \rho(\gamma)v=v. $$
Denote by $\mathcal{Z}_\lambda$ the set
$$ \mathcal{Z}_\lambda:=\{ g \in \G\ :\ \lambda(g)v=v\}.$$
Since $\lambda$ is irreducible, $\mathcal{Z}_\lambda \subsetneq \G$ and it is a real algebraic subset of $\G$. Indeed all irreducible representations can be included \footnote{This is a folklore theorem in the representation theory of compact Lie groups and can be proved in multiple ways see for example the discussion at {\it https://mathoverflow.net/questions/58633/does-every-irreducible-representation-of-a-compact-group-occur-in-tensor-product} and references given here.} in
$$\mathcal{I}^{\otimes k} \otimes     \widehat{\mathcal{I}}^{\otimes l}$$
for some $k,l\in \N$, and thus matrix coefficients are Laurent polynomials in the matrix coefficients of elements of $\G$. This translates (by taking real and imaginary parts) into real algebraic equations.
Since $\rho(\Gamma)$ is assumed to be Zariski dense, we must have
$$\rho(\Gamma)\cap \mathcal{Z}_\lambda^c\neq \emptyset,$$
and thus $\lambda\circ \rho$ cannot contain the trivial representation.
The sum over irreducible $\lambda$ being finite, we can apply the prime orbit theorem with remainder and Proposition \ref{prime2} to conclude the proof. $\square$

We point out that the leading term in the above asymptotic should follow directly from a dynamical result of Parry and Pollicott \cite{PP1}, provided that one can check the weak-mixing property
of the associated extension of the geodesic flow via the Zariski density assumption. We have chosen here a more self-contained and scenic road which gives a stronger result which might be useful in future applications.

\section{Proof of the main theorem}
\subsection{An equidistribution result}
In core of the proof of the main theorem, we will need to use the following equidistribution result which we state and prove here.
\begin{lem}
\label{equi1}
Let $0<\ell_0< \ell_1<\ldots <\ell_n$ be a sequence such that $\ell_n\rightarrow +\infty$ as $n\rightarrow +\infty$ and let $(a_n)_{n\geq 0}$ be a sequence of positive real numbers.
Consider the counting function
$$\Pi(x):=\sum_{\ell_n \leq x} a_n, $$
and assume that as $x\rightarrow +\infty$, we have $\Pi(x)\sim \kappa \frac{e^x}{x}$, for some $\kappa>0$. Assume that $g:\R^+\rightarrow \R^+$ is a continuous function such that $g(x)=\frac{x^2}{e^x}\left ( 1+O(1/x) \right)$ as $x\rightarrow +\infty$. Pick 
$\widehat{\psi}$ a positive valued, compactly supported $C^\infty$ function on $\R^+$. Fix $\alpha \in \R$. Then as $T\rightarrow +\infty$, we have
$$\sum_{\ell_n} a_n \cos^2(\alpha \ell_n) g(\ell_n)\widehat{\psi}\left(\frac{\ell_n}{T}\right)\sim \frac{\kappa T^2}{2} \int_0^\infty u \widehat{\psi}(u)du,$$
if $\alpha\neq 0$, while
$$\sum_{\ell_n} a_n \cos^2(\alpha \ell_n) g(\ell_n)\widehat{\psi}\left(\frac{\ell_n}{T}\right)\sim \kappa T^2 \int_0^\infty u \widehat{\psi}(u)du $$
if $\alpha=0$.
\end{lem}
\noindent {\it Proof}. First notice that for some constant $C>0$, we have the following bound valid for all $x\geq 0$,
$$\Pi(x)\leq \left \{ C\frac{e^x}{x}\ \mathrm{if}\ x\geq 1 \atop C\ \mathrm{if}\ 0\leq x\leq 1. \right. $$
Set 
$$S_T:= \sum_{\ell_n} a_n \cos^2(\alpha \ell_n)g(\ell_n)\widehat{\psi}\left(\frac{\ell_n}{T}\right).$$
By the asymptotic behaviour of $g$, we have as $T\rightarrow \infty$, 
$$S_T=\int_0^\infty \cos^2(\alpha x)\frac{x^2}{e^x}\widehat{\psi}\left(\frac{x}{T}\right)d\Pi(x) +O\left ( \int_0^{C_1 T} \frac{x}{e^x}d\Pi(x)   \right),$$
for some $C_1>0$ which depends on the support of $\widehat{\psi}$.
Using the above bound for $\Pi(x)$ and a summation by parts we get
$$  \int_0^{C_1 T} \frac{x}{e^x}d\Pi(x)= \int_0^{C_1 T} e^{-x}(x-1)\Pi(x)dx=O(T).$$
We are left to analyze 
$$S_T^0:= \int_0^\infty \cos^2(\alpha x)\frac{x^2}{e^x}\widehat{\psi}\left(\frac{x}{T}\right)d\Pi(x).$$
Using the identity
$$\cos^2(z)=\frac{\cos(2z)+1}{2},$$
we write
$$ S_T^0:=\widetilde{S}_T+S_T^1,$$
where 
$$ \widetilde{S}_T= \frac{1}{2}\int_0^\infty \frac{x^2}{e^x}\widehat{\psi}\left(\frac{x}{T}\right)d\Pi(x),\ S_T^1=\frac{1}{2}\int_0^\infty \cos(2\alpha x)\frac{x^2}{e^x}\widehat{\psi}\left(\frac{x}{T}\right)d\Pi(x).$$
We first deal with $\widetilde{S}_T$. A summation by parts yields
$$ \widetilde{S}_T=\frac{1}{2}\int_0^\infty e^{-x}(x^2-2x)\widehat{\psi}\left(\frac{x}{T}\right)\Pi(x)dx-\frac{2}{T}\int_0^\infty \frac{x^2}{e^x}\widehat{\psi}'\left(\frac{x}{T}\right)\Pi(x)dx.$$
We then observe that
$$\frac{1}{T}\int_0^\infty \frac{x^2}{e^x}\widehat{\psi}'\left(\frac{x}{T}\right)\Pi(x)dx=O\left( \frac{1}{T}\int_0^\infty x \left \vert\widehat{\psi}'\left(\frac{x}{T}\right) \right\vert dx   \right)=O(T),$$
while 
$$\int_0^\infty e^{-x}(x^2-2x)\widehat{\psi}\left(\frac{x}{T}\right)\Pi(x)dx=\int_0^\infty x^2e^{-x}\widehat{\psi}\left(\frac{x}{T}\right)\Pi(x)dx+O(T).$$
Therefore as $T\rightarrow +\infty$,
$$ \widetilde{S}_T=\frac{1}{2}\int_0^\infty x^2e^{-x}\widehat{\psi}\left(\frac{x}{T}\right)\Pi(x)dx+O(T). $$
By a change of variable,
we have
$$\int_0^\infty x^2e^{-x}\widehat{\psi}\left(\frac{x}{T}\right)\Pi(x)dx=T^2\int_0^\infty f_T(u)du, $$
where $f_T(u)=Tu^2 e^{-Tu}\widehat{\psi}(u) \Pi(Tu)$.
Notice that we have for fixed $u\geq 0$, 
$$\lim_{T\rightarrow \infty} f_T(u)= \kappa u\widehat{\psi}(u).$$
On the other hand we have for all $T$ large, for all $u\geq 0$, $$0\leq f_T(u)\leq Cu\widehat{\psi}(u),$$ therefore Lebesgue's dominated convergence theorem yields
$$\lim_{T\rightarrow \infty} \frac{ \widetilde{S}_T}{T^2}=\frac{\kappa}{2}\int_0^\infty u\widehat{\psi}(u)du,$$
which is the asymptotic formula expected. If $\alpha=0$, then the proof is done since $S_T^1= \widetilde{S}_T $. We now deal with $S_T^1$ assuming $\alpha\neq 0$.
A summation by parts combined with the crude bound on $\Pi(x)$ as above yields for large $T$
$$S_T^1=\frac{1}{2}\int_0^\infty \cos(2\alpha x)x^2e^{-x}\widehat{\psi}\left(\frac{x}{T}\right)\Pi(x)dx $$
$$+\alpha \int_0^\infty \sin(2\alpha x)x^2e^{-x}\widehat{\psi}\left(\frac{x}{T}\right)\Pi(x)dx+O(T).$$
Let us analyze the term
$$I_T:=\int_0^\infty \cos(2\alpha x)x^2e^{-x}\widehat{\psi}\left(\frac{x}{T}\right)\Pi(x)dx,$$
the other term involving the sine being similar. By a change of variable we get
$$I_T=T^3\int_0^\infty u^2 e^{-Tu}\cos(2\alpha Tu)\Pi(Tu)\widehat{\psi}(u)du$$
$$=T^2\int_0^\infty f_T(u)\cos(2\alpha Tu)du,$$
where $f_T(u)=Tu^2e^{-Tu}\Pi(Tu)\widehat{\psi}(u)$ as above.
We then write
$$\frac{I_T}{T^2}=\int_0^\infty (f_T(u)-\kappa u\widehat{\psi}(u))\cos(2\alpha Tu)du+\kappa\int_0^\infty u\widehat{\psi}(u) \cos(2 \alpha Tu)du,$$
and a single integration by parts shows that
$$\int_0^\infty u\widehat{\psi}(u) \cos(2 \alpha Tu)du=O\left( \frac{1}{T} \right),$$
while
$$\left \vert \int_0^\infty (f_T(u)-\kappa u\widehat{\psi}(u))\cos(2\alpha Tu)du\right \vert \leq \int_{\R^+}\vert f_T(u)-\kappa u\widehat{\psi}(u)\vert du.$$
By applying Lebesgue's theorem as before, we finally conclude that
$$\lim_{T\rightarrow +\infty} \frac{I_T}{T^2}=0, $$
and thus
$$ \lim_{T\rightarrow +\infty} \frac{S_T}{T^2}=\lim_{T\rightarrow +\infty} \frac{\widetilde{S}_T}{T^2}=\frac{\kappa}{2}\int_0^\infty u\widehat{\psi}(u)du,$$
and the proof is done. $\square$

The rest of this section is devoted to the proof of Theroem \ref{main} and is naturally subdivided into three main steps. We recall that we investigate the spectrum of
$\Delta_{n,\rho}$ where $\rho:\Gamma \rightarrow GL(V_\rho)$ is a finite dimensional unitary representation. For simplicity we will denote the character of $\rho$ by
$$\chi(\gamma):=\mathrm{Tr}(\rho(\gamma)).$$

Let $\psi$ be a $C^\infty$ function whose Fourier transform
$\widehat{\psi}=\widehat{\psi}$ is in $C_0^\infty(\R)$, is even, and takes its values in $\R^+$. Notice that $\psi$ is given by
$$\psi(r)=\int_\R \widehat{\psi}(u)e^{iru}du, $$
and is therefore even too.
We will apply the (twisted) trace formula to $$h(r)=\psi(L(r-\alpha))+\psi(L(r+\alpha)),$$ where 
$\alpha$ is fixed and $L$ will be taken large. We then have
$$\widehat{h}(u)= \frac{2\cos(\alpha u)}{L}\widehat{\psi}\left( \frac{u}{L} \right),$$
so that 
$$N_n(L)=(g-1)n\mathrm{dim}(V_\rho)\int_\R  \left(\psi(L(r-\alpha))+\psi(L(r+\alpha)) \right)r\tanh(\pi r)dr +N_n^0(L),$$
where we have set
$$ N_n^0(L)=\frac{2}{L}\sum_{\gamma \in \mathcal{P}} \sum_{k\geq 1}  \chi(\gamma^k)F_n(\gamma^k) \frac{\ell(\gamma)\widehat{\psi}(k\ell(\gamma)/L)}{2\sinh(k\ell(\gamma)/2)}\cos(\alpha k \ell(\gamma)).$$
Notice that if $\gamma$ is a representative of a primitive conjugacy class, then $\gamma^{-1}$ represents another conjugacy class which corresponds to a closed geodesic with reversed orientation. Therefore the character $\chi$ breaks the {\it time reversal symmetry} if and only if it can distinguish the orientation of some closed geodesics. As emphasized by the physical intuition,
it will prove more convenient to re-write the above sum over {\it non-oriented} primitive closed geodesics: we set $\mathcal{P}_0=\mathcal{P}/\{\mathrm{Id}, \mathrm{Inv} \}$,
where $\mathrm{Inv}:\mathcal{P}\rightarrow \mathcal{P}$ is the involution given by the inversion, acting on conjugacy classes. We then have
$$N_n^0(L)=\frac{2}{L}\sum_{\gamma \in \mathcal{P}_0} \sum_{k\geq 1}  \left(\chi(\gamma^k)+\overline{\chi(\gamma^k)}\right) F_n(\gamma^k) \frac{\ell(\gamma)\widehat{\psi}(k\ell(\gamma)/L)}{2\sinh(k\ell(\gamma)/2)}\cos(\alpha k \ell(\gamma)).$$

\subsection{Expectation and variance in the large $n$ limit}
We first observe that by taking expectations and using the fact that for all $k$,  
$$\lim_{n\rightarrow } \E_n(F_n(\gamma^k))=d(k)$$ exists, we obtain that for large $n$ (here $\alpha,L$ are fixed),
$$\E_n(N_n(L))=(g-1)n\mathrm{dim}(V_\rho)\int_\R  \left(\psi(L(r-\alpha))+\psi(L(r+\alpha)) \right)r\tanh(\pi r)dr+O_L(1),$$
which is a smooth probabilistic Weyl law in the large $n$ regime. The quantum variance is 
$$\mathbb{V}_n(L):=\mathbb{E}_n \left (\left \vert N_n(L)-\mathbb{E}_n(N_n(L)) \right \vert^2\right)=\mathbb{E}_n \left (\left \vert N_n^0(L)-\mathbb{E}_n(N_n^0(L)) \right \vert^2\right).$$
for simplicity set 
$$G_k(\gamma):=\left(\chi(\gamma^k)+\overline{\chi(\gamma^k)}\right)\cos(\alpha k \ell(\gamma))\frac{\ell(\gamma)\widehat{\psi}(k\ell(\gamma)/L)}{2\sinh(k\ell(\gamma)/2)} , $$
so that we have
$$\mathbb{V}_n(L)=$$
$$\frac{4}{L^2} \sum_{\gamma_1,k_1} \sum_{\gamma_2,k_2} G_{k_1}(\gamma_1)\overline{G_{k_2}(\gamma_2)}
\E_n\left [ \left(F_n(\gamma_1^{k_1})-\E_n(F_n(\gamma_1^{k_1}))\right) \left(F_n(\gamma_2^{k_2})-\E_n(F_n(\gamma_2^{k_2}))\right)\right].$$
The parameter $L$ being fixed, we observe by Proposition \ref{Stat1} $(1)+(2)$ that whenever $\gamma_1\neq \gamma_2 \in \mathcal{P}_0$ we have
$$\lim_{n\rightarrow +\infty} \E_n\left [ \left(F_n(\gamma_1^{k_1})-\E_n(F_n(\gamma_1^{k_1}))\right) \left(F_n(\gamma_2^{k_2})-\E_n(F_n(\gamma_2^{k_2}))\right)\right]=0,$$
therefore by Proposition \ref{Stat1} $(3)$ we get
$$\lim_{n\rightarrow +\infty} \mathbb{V}_n(L):=\mathcal{S}(L):=\frac{4}{L^2} \sum_{k_1,k_2} \mathcal{V}(k_1,k_2)S_{k_1,k_2}(L),$$
where 
$$ S_{k_1,k_2}(L):=\sum_{\gamma \in \mathcal{P}_0} G_{k_1}(\gamma)\overline{G_{k_2}(\gamma)}.$$
We have obtained a "near diagonal" sum, and the goal of the remaining subsections is to analyze the asymptotic behaviour as $L\rightarrow +\infty$ of each term $S_{k_1,k_2}(L)$.
\subsection{The non-primitive contribution}
The non-primitive contribution corresponds to the sum over integers $k_1,k_2$ such that $k_1+k_3\geq 3$. We warn the reader that in the estimates below, 
the implied constants do not depend on $L,k_1,k_2$ and may change from line to line. We will therefore use for simplicity Vinogradov's notation $f(T)\ll g(T)$ which means that $f(T)=O(g(T))$ as $T\rightarrow +\infty$.
By a crude bound we have 
$$\vert S_{k_1,k_2}(L)\vert\ll \sum_{\gamma \in \mathcal{P}_0} e^{-\frac{k_1+k_2}{2}\ell(\gamma)}\ell(\gamma)^2.$$
so that it is actually uniformly bounded with respect to $L$.
 By a Sieltjes integration by parts we have
$$\vert S_{k_1,k_2}(L)\vert\ll\int_0^\infty  e^{-\frac{k_1+k_2}{2}u}\left ( \frac{k_1+k_2}{2}u^2-2u\right)N^0(u)du$$
$$\ll(k_1+k_2)\int_0^\infty  e^{-(\frac{k_1+k_2}{2}-1)u}u^2du,$$
by using a crude bound on $N^0(u)\leq C e^u$.
Using the fact that for all $\alpha>0$ we have
$$\int_0^\infty e^{-\alpha t}t^2dt=\frac{1}{\alpha^3},$$
we get for all $k_1+k_2\geq 3$,
$$\vert S_{k_1,k_2}(L)\vert\ll \frac{1}{(k_1+k_2)^2}.$$
We are now ready to estimate
$$S_{NP}(L):=\sum_{k_1+k_2\geq 3} \mathcal{V}(k_1,k_2)  S_{k_1,k_2}(L).$$
We will show that $S_{NP}(L)\ll \log^2(L)$, following a computation that was suggested to us by Zeev Rudnick.
First observe that since $\widehat{\psi}$ has compact support, the sum over $k_1,k_2$ is limited to the range
$$k_1,k_2\ll L.$$
 Notice that we have also
$$\mathcal{V}(k_1,k_2)\ll \sigma( \gcd(k_1,k_2) ),$$
where
$$\sigma(n)=\sum_{d\vert n} d $$
is the sum of divisors function. Hence we have
$$ S_{NP}(L)\ll \sum_{k_1,k_2 \ll L} \frac{\sigma( \gcd(k_1,k_2) )}{(k_1+k_2)^2}.$$
By writing for $j=1,2$, $k_j=Dm_j$ with $\gcd(m_1,m_2)=1$, we get
$$S_{NP}(L)\ll \sum_{D}\sum_{m_1D,m_2D\ll L}  \frac{\sigma( D)}{D^2(m_1+m_2)^2}$$
$$\ll \left( \sum_{D\ll L}  \frac{\sigma(D)}{D^2}\right) \left( \sum_{m_1,m_2\ll L}  \frac{1}{(m_1+m_2)^2}\right).$$
It is easy to check that
$$\left( \sum_{m_1,m_2\ll L}  \frac{1}{(m_1+m_2)^2}\right)\ll \log(L).$$
On the other hand, we can write
$$\sigma(D)=\sum_{ ab=D} a$$
so that
$$ \sum_{D\ll L}  \frac{\sigma(D)}{D^2}\ll \sum_{a,b\ll L} \frac{a}{(ab)^2}\ll \sum_{a\ll L} \frac{1}{a} \sum_{b=1}^\infty \frac{1}{b^2}\ll \log(L).$$

In a nutshell, we have reached 
$$\lim_{n\rightarrow +\infty} \mathbb{V}_n(L)=\mathcal{S}(L)=\frac{4}{L^2} \sum_{k_1,k_2} \mathcal{V}(k_1,k_2)S_{k_1,k_2}(L) $$
$$=\frac{4}{L^2}S_{1,1}(L)+O\left( \frac{(\log(L))^2}{L^{2} }\right),$$
where
$$S_{1,1}(L)=\sum_{\gamma \in \mathcal{P}_0}   \left(\chi(\gamma)+\overline{\chi(\gamma)}\right)^2 \frac{\ell(\gamma)^2\widehat{\psi}^2(\ell(\gamma)/L)}{4\sinh^2(\ell(\gamma)/2)}\cos^2(\alpha \ell(\gamma)).$$
\subsection{Applying the equidistribution Lemma}
We now know that
$$\lim_{L\rightarrow +\infty}\lim_{n\rightarrow +\infty} \mathbb{V}_n(L)=\lim_{L\rightarrow +\infty}  \frac{4}{L^2}S_{1,1}(L),$$
thus we need to determine the asymptotic behaviour when $L\rightarrow +\infty$ of $S_{1,1}(L)$. We assume that $\alpha \neq 0$.
We will therefore apply Lemma \ref{equi1} with
$$\Pi(x)=\sum_{\gamma \in \mathcal{P}_0 \atop \ell(\gamma)\leq x}  \left(\chi(\gamma)+\overline{\chi(\gamma)}\right)^2,\ g(x)=\frac{x^2}{4\sinh^2(x/2)}.$$
Writing
$$ \Pi(x)=\sum_{\gamma \in \mathcal{P}_0 \atop \ell(\gamma)\leq x} \chi^2(\gamma)+\sum_{\gamma \in \mathcal{P}_0 \atop \ell(\gamma)\leq x} \overline{\chi}^2(\gamma)
+\sum_{\gamma \in \mathcal{P}_0 \atop \ell(\gamma)\leq x} 2,$$
we observe that two cases can occur. Either $\chi^2\equiv1$ and we have by the prime orbit theorem (recall that we are summing over $\mathcal{P}_0$)
$$\Pi(x)\sim 2\frac{e^x}{x},$$
which by applying Lemma \ref{equi1} gives
$$\lim_{L\rightarrow +\infty}  \frac{4}{L^2}S_{1,1}(L)=4\int_0^\infty u \widehat{\psi}^2(u)du=2\int_\R \vert u \vert \left ( \widehat{\psi}(u))\right)^2 du=\Sigma^2_{\mathrm{GOE}}(\psi). $$
Or we have $\chi^2\not \equiv 1$ and we can apply Proposition \ref{prime2} to get as $x\rightarrow +\infty$
$$\Pi(x)=\frac{e^x}{x}+ O\left ( \frac{e^x}{x^2} \right )+O\left( e^{\beta_0 x}\right),$$
for some $\beta_0<1$. Therefore as $x\rightarrow +\infty$ we have
$$\Pi(x)\sim \frac{e^x}{x},$$
and applying Lemma \ref{equi1} gives now
$$\lim_{L\rightarrow +\infty}  \frac{4}{L^2}S_{1,1}(L)=2\int_0^\infty u \widehat{\psi}^2(u)du=\int_\R \vert u \vert \left ( \widehat{\psi}(u))\right)^2 du=\Sigma^2_{\mathrm{GUE}}(\psi). $$
Theorem \ref{main1} is proved. We now show how Theorem \ref{main2} follows. In view of the above dicussion, all we have to do is to compute the asymptotics of
$$\Pi(x)=\sum_{\gamma \in \mathcal{P}_0 \atop \ell(\gamma)\leq x}  \left(\mathrm{Tr}(\rho(\gamma))+\overline{\mathrm{Tr}(\rho(\gamma))}\right)^2.$$
We can apply Proposition \ref{equi2} with $\G$ as in Theorem \ref{main2} and 
$$f(g)=\left(\mathrm{Tr}(g)+\overline{\mathrm{Tr}(g)}\right)^2$$ which says that
$$\Pi(x)\sim  \frac{1}{2}\left (\int_{\G}f(g)dg\right)\frac{e^x}{x}.$$
Assume that $\G \subset U(N)$. First we can observe that
$$\int_\G f(g)dg=\int_\G (\mathrm{Tr}(g))^2dg+\overline{\int_\G (\mathrm{Tr}(g))^2dg } +2 \int_\G \vert \mathrm{Tr}(g)\vert^2dg.$$
As done previously in $\S 4$, let $\mathcal{I}:\G \rightarrow GL(\C^N)$ denote the standard representation of $\G$ acting on $\C^N$ by left matrix multiplication.
Since $\G$ is not a direct product, $\mathcal{I}$ is actually irreducible, hence Schur orthogonality says that
$$\int_\G  \vert \mathrm{Tr}(g)\vert^2dg=\int_\G \vert \mathrm{Tr}(\mathcal{I}(g))\vert^2dg=1.$$
On the other hand,
$$\int_\G (\mathrm{Tr}(g))^2dg=\int_\G \mathrm{Tr}(\mathcal{I}\otimes \mathcal{I}(g))dg=d_{\mathrm{trivial}},$$
where $d_{\mathrm{trivial}}$ denotes the multiplicity of the trivial representation in the irreducible decomposition of $\mathcal{I}\otimes \mathcal{I}$, in particular it is a non-negative integer.  Since
$$\left \vert  \int_\G (\mathrm{Tr}(g))^2dg\right \vert \leq \int_\G  \vert \mathrm{Tr}(g)\vert^2dg=1,$$
we must have $d_{\mathrm{trivial}}\in \{ 0,1\}$. To summarize,
$$\int_\G f(g)dg=2\int_\G (\mathrm{Tr}(g))^2dg+ 2\in \{2,4\},$$
which shows that the only possible limits of the spectral variance are $\Sigma^2_{\mathrm{GOE}}$ and $\Sigma^2_{\mathrm{GUE}}$. Since
$$\int_\G (\mathrm{Tr}(g))^2dg \in \R, $$
we have
$$0\leq \int_\G (\mathrm{Tr}(g))^2dg=\int_\G (\Re(\mathrm{Tr}(g)))^2dg-\int_\G (\Im(\mathrm{Tr}(g)))^2dg\leq \int_\G  \vert \mathrm{Tr}(g)\vert^2dg=1.$$
Clearly if $g\mapsto \mathrm{Tr}(g)\in \R$, then
$$\int_\G (\mathrm{Tr}(g))^2dg=\int_\G  \vert \mathrm{Tr}(g)\vert^2dg=1,$$
and we are in the GOE case. If there exists $g\in \G$ such that $\Im(\mathrm{Tr}(g))\neq 0$, then a simple continuity argument shows that
$$\int_\G (\mathrm{Tr}(g))^2dg< \int_\G  \vert \mathrm{Tr}(g)\vert^2dg=1,$$
and therefore
$$\int_\G (\mathrm{Tr}(g))^2dg=0, $$
we are in the GUE case. Applying this criterion to the classical groups shows that $\G=SO(N),Sp(N)$ will fall in the GOE case since traces are always real.
On the other hand, for $\G=U(N)$ and $\G=SU(N)$ with $N\geq 3$, traces can take complex values and we get the GUE case. Notice that $Sp(1)=SU(2)$.
\bigskip
Theorem \ref{main2} is now proved. $\square$.

The reader might wonder if representations $\rho:\Gamma \rightarrow \G$ with Zariski dense images actually exist. A paper of Kishore \cite{kishore}, using results of Breuillard, Green, Gurnalik and Tao, shows that if $\G$ is real algebraic and {\it semi-simple}, then (Theorem 2.7), such a homomorphism always exists. This is enough to deal with $\G=SU(N),SO(N),SP(N)$. More Generally, in \cite{KP}, they show (Corollary 2.4) using a dimensional argument, that for any real algebraic {\it reductive} Lie Group $\G$, the set of representations with Zariski dense image is non-empty (and actually dense) provided that
$$g\geq \left ( \mathrm{dim}_\R(\G)\right)^2,$$
where $g$ is the genus of $\Gamma \backslash \H$ and $\mathrm{dim}_\R$ is the real dimension of $\G$. This allows to build representations with zariski dense images in $U(N)$
provided $g\geq N^4$.


\bigskip
\begin{thebibliography}{10}
 \bibitem{AM}
Nalini Anantharaman and Laura Monk.
\newblock A high-genus asymptotic expansion of {W}eil-{P}etersson volume
  polynomials.
\newblock {\em J. Math. Phys.}, 63(4):Paper No. 043502, 26, 2022.

\bibitem{Berry1}
M.~V. Berry.
\newblock Fluctuations in numbers of energy levels.
\newblock In {\em Stochastic processes in classical and quantum systems
  ({A}scona, 1985)}, volume 262 of {\em Lecture Notes in Phys.}, pages 47--53.
  Springer, Berlin, 1986.

\bibitem{Berry2}
Michael Berry.
\newblock Semiclassical mechanics of regular and irregular motion.
\newblock In {\em Chaotic behavior of deterministic systems ({L}es {H}ouches,
  1981)}, pages 171--271. North-Holland, Amsterdam, 1983.

\bibitem{BS}
E.~B. Bogomolny, B.~Georgeot, M.-J. Giannoni, and C.~Schmit.
\newblock Chaotic billiards generated by arithmetic groups.
\newblock {\em Phys. Rev. Lett.}, 69(10):1477--1480, 1992.

\bibitem{Bump}
Daniel Bump.
\newblock {\em Lie groups}, volume 225 of {\em Graduate Texts in Mathematics}.
\newblock Springer-Verlag, New York, 2004.

\bibitem{Chevalley}
Claude Chevalley.
\newblock {\em Theory of {L}ie {G}roups. {I}}.
\newblock Princeton Mathematical Series, vol. 8. Princeton University Press,
  Princeton, N. J., 1946.

\bibitem{DM1}
Freeman~J. Dyson and Madan~Lal Mehta.
\newblock Statistical theory of the energy levels of complex systems. {IV}.
\newblock {\em J. Mathematical Phys.}, 4:701--712, 1963.

\bibitem{Hejhal1}
Dennis~A. Hejhal.
\newblock {\em The {S}elberg trace formula for {${\rm PSL}(2,R)$}. {V}ol. {I}}.
\newblock Lecture Notes in Mathematics, Vol. 548. Springer-Verlag, Berlin-New
  York, 1976.

\bibitem{Huber}
Heinz Huber.
\newblock Zur analytischen {T}heorie hyperbolischen {R}aumformen und
  {B}ewegungsgruppen.
\newblock {\em Math. Ann.}, 138:1--26, 1959.

\bibitem{KP}
Inkang Kim and Pierre Pansu.
\newblock Density of {Z}ariski density for surface groups.
\newblock {\em Duke Math. J.}, 163(9):1737--1794, 2014.

\bibitem{kishore}
Krishna Kishore.
\newblock Representation variety of surface groups.
\newblock {\em Proc. Amer. Math. Soc.}, 146(3):953--959, 2018.

\bibitem{LS1}
Martin~W. Liebeck and Aner Shalev.
\newblock Fuchsian groups, coverings of {R}iemann surfaces, subgroup growth,
  random quotients and random walks.
\newblock {\em J. Algebra}, 276(2):552--601, 2004.

\bibitem{LuoSarnak}
W.~Luo and P.~Sarnak.
\newblock Number variance for arithmetic hyperbolic surfaces.
\newblock {\em Comm. Math. Phys.}, 161(2):419--432, 1994.

\bibitem{MNP}
Michael Magee, Fr{\'e}d{\'e}ric Naud, and Doron Puder.
\newblock A random cover of a compact hyperbolic surface has relative spectral
  gap {$\frac{3}{16}-\epsilon$}.
\newblock {\em To appear in GAFA}, 2022.

\bibitem{MP1}
Michael Magee and Doron Puder.
\newblock The asymptotic statistics of random covering surfaces.
\newblock {\em Preprint}, 2020.

\bibitem{MehtaBook}
Madan~Lal Mehta.
\newblock {\em Random matrices}, volume 142 of {\em Pure and Applied
  Mathematics (Amsterdam)}.
\newblock Elsevier/Academic Press, Amsterdam, third edition, 2004.

\bibitem{Mirza}
Maryam Mirzakhani.
\newblock Simple geodesics and {W}eil-{P}etersson volumes of moduli spaces of
  bordered {R}iemann surfaces.
\newblock {\em Invent. Math.}, 167(1):179--222, 2007.

\bibitem{MirPet}
Maryam Mirzakhani and Bram Petri.
\newblock Lengths of closed geodesics on random surfaces of large genus.
\newblock {\em Comment. Math. Helv.}, 94(4):869--889, 2019.

\bibitem{Nica}
Alexandru Nica.
\newblock On the number of cycles of given length of a free word in several
  random permutations.
\newblock {\em Random Structures Algorithms}, 5(5):703--730, 1994.

\bibitem{PP1}
William Parry and Mark Pollicott.
\newblock The {C}hebotarov theorem for {G}alois coverings of {A}xiom {A} flows.
\newblock {\em Ergodic Theory Dynam. Systems}, 6(1):133--148, 1986.

\bibitem{PS1}
Ralph Phillips and Peter Sarnak.
\newblock Geodesics in homology classes.
\newblock {\em Duke Math. J.}, 55(2):287--297, 1987.

\bibitem{PZ}
Doron Puder and Tomer Zimhoni.
\newblock Local statistics of random permutations from free products.
\newblock {\em Preprint}, 2022.

\bibitem{Rudnick}
Zeev Rudnick.
\newblock Spectral statistics on the moduli space of surfaces of large genus.
\newblock {\em Preprint}, 2022.

\bibitem{Sarnak1}
Peter~Clive Sarnak.
\newblock {\em P{RIME} {GEODESIC} {THEOREMS}}.
\newblock ProQuest LLC, Ann Arbor, MI, 1980.
\newblock Thesis (Ph.D.)--Stanford University.

\bibitem{Sunada}
Toshikazu Sunada.
\newblock Unitary representations of fundamental groups and the spectrum of
  twisted {L}aplacians.
\newblock {\em Topology}, 28(2):125--132, 1989.

\bibitem{Venkov}
Alexei~B. Venkov.
\newblock {\em Spectral theory of automorphic functions and its applications},
  volume~51 of {\em Mathematics and its Applications (Soviet Series)}.
\newblock Kluwer Academic Publishers Group, Dordrecht, 1990.
\newblock Translated from the Russian by N. B. Lebedinskaya.


\end{thebibliography}
\end{document}